\documentclass{article}

\usepackage{amsthm,amsmath}
\usepackage{graphicx}
\usepackage{enumerate}
\usepackage{amssymb} % allow blackboard bold (aka N,R,Q sets)

\RequirePackage[numbers]{natbib}
\RequirePackage[colorlinks,citecolor=blue,urlcolor=blue]{hyperref}

\def\Item$#1${\item {\hfill $\displaystyle#1$ \hfill\refstepcounter{equation}(\theequation)}}

\begin{document}

\title{Sums of Squares\\ Bijective Parameter Representation}
\author{Walter Wyss}
\date{}
\maketitle

\begin{abstract}
A pair of n-dimensional vectors x, y $\in$ $\mathbb{R}^n$ of equal length \\
\begin{equation*}
x=(x_1, x_2, \dotsb, x_n) , \hspace{0.5cm} y=(y_1, y_2, \dotsb, y_n) \\
\end{equation*}
\begin{equation}
x_1^2 + x_2^2 + \dotsb + x_n^2 = y_1^2 + y_2^2 + \dotsb + y_n^2 \\
\end{equation}
\\
generate a convenient parameterization that uniquely reproduces the vectors. If the components of the vectors are rational numbers, the parameters are rational and vice versa. Multiplying with the greatest common denominator we get the complete solution, up to scaling and absolute values, of equation (1) in positive integers. \\
We give examples where some of the components are equal or coincide. We also apply our representation to the parallelogram equation \\
\begin{equation}
2u_1^2 + 2u_2^2 = u_3^2 + u_4^2
\end{equation}
and to the equation
\begin{equation}
x^2 + y^2 + z^2 = 3w^2
\end{equation}
\end{abstract}

\section {Parameterization. \\} 
For x,y $\in$ $\mathbb{R}^n$  \\
\begin{equation*}
x=(x_1, x_2, \dotsb, x_n) ,  \hspace{0.5cm} y=(y_1, y_2, \dotsb, y_n) \\
\end{equation*}
\\
we have the scalar product  \\
\begin{equation}
(x,y) = \sum_{i=1}^{n} x_iy_i \\
\end{equation}
\\
and the norm $\Vert x \Vert$ through \\
\begin{equation}
\Vert x \Vert^2 = \sum_{i=1}^{n} x_1^2 \\
\end{equation}
\\
Equation (1) means that the vectors x and y have equal length \\

\begin{equation}
\Vert x \Vert = \Vert y \Vert \\
\end{equation}
\\
Equivalently to (1) we have \\

\begin{equation}
x_1^2 - y_1^2 + \dotsb + x_k^2 - y_k^2 + \dotsb + x_n^2 - y_n^2 = 0 \\
\end{equation}
\\
We now introduce the quantities $s_i$ and $d_i$ by \\

\begin{equation}
x_i + yi = 2s_i , \hspace{0.5cm} x_i - y_i = 2d_i, \hspace{0.5cm} i = 1, 2, \dotsb, n \\
\end{equation}
\\
giving \\

\begin{equation}
x_i = s_i + d_i , \hspace{0.5cm} y_i = s_i - d_i \\
\end{equation}
\\
In vector notation \\

\begin{equation}
s = (s_1, \dotsb, s_n) , \hspace{0.5cm} d = (d_1, \dotsb, d_n) \\
\end{equation}
\\
we have \\

\begin{equation}
x = s + d, \hspace{0.5cm} y = s - d  \\
\end{equation}
\\
From (7) we find \\

\begin{equation}
(s,d) = 0 \\
\end{equation}

\begin{equation}
\Vert x \Vert^2 = \Vert s \Vert^2 + \Vert d \Vert^2 \\
\end{equation}
\\
Thus the vectors s and d are orthogonal. \\
\\
Given the vector \\

\begin{equation}
s = (s_1, \dotsb, s_n) , \hspace{0.5cm} s_1 \neq 0 \\
\end{equation}
\\
we look at the vectors, $k = 2, 3, \dotsb, n$ \\

\begin{equation}
\begin{aligned}
e_2& = ( - s_2, s_1, 0, \dotsb, 0 ) \\
e_k& = ( - s_k, 0, \dotsb s_1, 0, \dotsb, 0), \hspace{0.5cm} s_1 \text{in the k\textsuperscript{th} place} \\
e_n& = ( - s_n, 0, \dotsb, 0, s_1) \\
\end{aligned}
\end{equation}
\\
These are (n-1) vectors that are linearly independent and orthogonal to the vector s. Thus they span the vector space orthogonal to the vector s. \\
Since the vector d is orthogonal to the vector s, it can be written as \\

\begin{equation}
d = \sum_{k=2}^{n} \lambda_k e_k  \\
\end{equation}
\\
The parameterization of (1) is now given by \\

\begin{equation}
x = s + \sum_{k=2}^{n} \lambda_k e_k  \\
\end{equation}

\begin{equation}
y = s - \sum_{k=2}^{n} \lambda_k e_k  \\
\end{equation}
\\
or in components \\

\begin{equation}
\begin{aligned}
x_1& = s_1 - \sum_{k=2}^{n} \lambda_k s_k  \\
x_k& = s_k + \lambda_k s_1 \hspace{0.5cm}, k=2, \dotsb, n \\
\end{aligned}
\end{equation}
\\

\begin{equation}
\begin{aligned}
y_1& = s_1 + \sum_{k=2}^{n} \lambda_k s_k  \\
y_k& = s_k - \lambda_k s_1 \hspace{0.5cm}, k=2, \dotsb, n \\
\end{aligned}
\end{equation}
\\
The (2n-1) parameters $s_1, \dotsb, s_n, \lambda_2, \dotsb, \lambda_n$ are now our convenient parameters giving us the components $x_i, y_i, i=1, \dotsb$ n. \\
Conversely, having the components $x_i$, $y_i$, we find the parameters by \\ 

\begin{equation}
2s_i = x_i + y_i,  \hspace{0.5cm} i=1, \dotsb, n \\
\end{equation}
 
\begin{equation}
2s_1\lambda_k = x_k - y_k,  \hspace{0.5cm} k=2, \dotsb, n \\
\end{equation}
\\
If the components $x_i, y_i$ are rational, the parameters $s_i, \lambda_k$ are rational and if the parameters $s_i, \lambda_k$ are rational the components $x_i, y_i$ are rational due to (19) and (20). There is a one-to-one correspondence between the parameters and the components.\\
\\

\section{Examples.}

\begin{enumerate}[a)]
 \Item  $ x_1^2 + x_2^2 = y_1^2 + y_2^2 $ \\
\\
	\begin{equation}
	\begin{aligned}
	x_1& = s_1 - \lambda_2s_2 \hspace{0.5cm} y_1 = s_1 + \lambda_2s_2 \\
	x_2& = s_2 + \lambda_2s_1 \hspace{0.5cm} y_2 = s_2 - \lambda_2s_1 \\
	\end{aligned}
	\end{equation}
\\
	Conversely \\
\\
	\begin{equation}
	\begin{aligned}
	2&s_1 = x_1 + y_1, \hspace{0.5cm} 2s_2 = x_2 + y_2 \\
	\end{aligned}
	\end{equation}

	\begin{equation}
	\begin{aligned}
	2&s_1 \lambda_2 = x_2 - y_2 \\
	\end{aligned}
	\end{equation}
\\
 \Item  $ x_1^2 + x_2^2= y_1^2 $ \\
\\
	Then $y_2=0$, resulting in $s_2=\lambda_2s_1$ and thus \\
	\begin{equation}
	\begin{aligned}
	x_1& = s_1 [1 - \lambda_2^2], \hspace{0.5cm} y_1 = s_1 [1+ \lambda_2^2] \\
	\end{aligned}
	\end{equation}

	\begin{equation}
	\begin{aligned}
	x_2& = 2s_1 \lambda_2 \\
	\end{aligned}
	\end{equation}
\\
	Conversely \\
\\
	\begin{equation}
	\begin{aligned}
	2&s_1 = x_1 + y_1, \hspace{0.5cm} 2s_1 \lambda_2 = x_2 \\
	\end{aligned}
	\end{equation}
\\
 \Item  $ x_1^2 + x_2^2 + \dotsb + x_n^2 = y_1^2 $ \\
\\
	Then $y_k = 0, k=2, \dotsb, n$ , resulting in $s_k=\lambda_ks_1$ \\
	and thus \\
	\begin{equation}
	\begin{aligned}
	x_1& = s_1 [ 1 - \lambda_2^2 - \dotsb - \lambda_n^2], \hspace{0.5cm} y_1 = s_1 [ 1 + \lambda_2^2 + \dotsb + \lambda_n^2] \\
	x_k& = 2s_1\lambda_k, \hspace{0.5cm} k = 2, \dotsb, n \\
	\end{aligned}
	\end{equation}
\\
	Conversely \\
\\
	\begin{equation}
	\begin{aligned}
	2&s_1 = x_1 + y_1, \hspace{0.5cm} 2s_1\lambda_k = x_k, \hspace{0.5cm} k=2, \dotsb, n \\
	\end{aligned}
	\end{equation}
\\

\end{enumerate}

\section{Applications.}

\begin{enumerate}[a)]
 \item Parallelogram equation [1] \\
\\
	\begin{equation}
	\begin{aligned}
	2u_1^2 + 2u_2^2 = u_3^2 + u_4^2 \\
	\end{aligned}
	\end{equation}
\\
We could use the parameterization (19) (20) for n=4, but it is more convenient to introduce \\
\\
	\begin{equation}
	\begin{aligned}
	2u_+ = u_4 + u_3, \hspace{0.5cm} 2u_- = u_4 - u_3 \\
	\end{aligned}
	\end{equation}
\\
Then\\
	\begin{equation}
	\begin{aligned}
	u_3 = u_+ - u_-, \hspace{0.5cm} u_4 = u_+ + u_- \\
	\end{aligned}
	\end{equation}
\\
and \\
\\
	\begin{equation}
	\begin{aligned}
	u_1^2 + u_2^2 = u_+^2 + u_-^2 \\
	\end{aligned}
	\end{equation}
\\
We now use (23) and get the bijective parameter representation \\
\\
	\begin{equation}
	\begin{aligned}
	u_1& = s_1 - \lambda_2s_2, \hspace{0.5cm} u_+ = s_1 + \lambda_2s_2 \\
	\end{aligned}
	\end{equation}

	\begin{equation}
	\begin{aligned}
	u_2& = s_2 + \lambda_2s_1, \hspace{0.5cm} u_- = s_2 - \lambda_2s_1 \\
	\end{aligned}
	\end{equation}
\\
or
\\
	\begin{equation}
	\begin{aligned}
	u_3& = s_1 + \lambda_2s_2 - s_2 + \lambda_2s_1 \\
	\end{aligned}
	\end{equation}

	\begin{equation}
	\begin{aligned}
	u_4& = s_1 + \lambda_2s_2 + s_2 - \lambda_2s_1 \\
	\end{aligned}
	\end{equation}
\\
Conversely, from (24) \\
\\
	\begin{equation}
	\begin{aligned}
	2s_1 = u_1 + u_+, \hspace{0.5cm} 2s_2 = u_2 +u_-, \hspace{0.5cm} 2s_1\lambda_2 = u_2 - u_- \\
	\end{aligned}
	\end{equation}
\\
or
\\
	\begin{equation}
	\begin{aligned}
	4s_1 = 2u_1 + u_4 + u_3, \hspace{0.5cm} 4s_2 = 2u_2 + u_4 - u_3 \\
	\end{aligned}
	\end{equation}

	\begin{equation}
	\begin{aligned}
	4s_1\lambda_2 = 2u_2 - u_4 + u_3 \\
	\end{aligned}
	\end{equation}
\\
Now introducing m, n, u through \\
\\
	\begin{equation}
	\begin{aligned}
	s_1 = u, \hspace{0.5cm} s_2 = nu, \hspace{0.5cm} \lambda_2 = m \\
	\end{aligned}
	\end{equation}

we find the bijective parameter representation of the parallelogram equation (34) in terms of the parameters m, n, u, as \\

	\begin{equation}
	\begin{aligned}
	u_1& = (1-mn)u  \\
	u_2& = (m+n)u  \\
	u_3& = (1+mn-n+m)u  \\
	u_4& = (1+mn+n-m)u  \\
	\end{aligned}
	\end{equation}
\\
Conversely\\
\\
	\begin{equation}
	\begin{aligned}
	4u& = 2u_1 + u_4 + u_3 \\
	m& = \frac{2u_2 - u_4 + u_3} {4u} \\
	n& = \frac{2u_2 + u_4 - u_3} {4u} \\
	\end{aligned}
	\end{equation}
\\
 \item The equation \\
\\
	\begin{equation}
	\begin{aligned}
	x^2& + y^2 + z^2 = 3w^2 \
	\end{aligned}
	\end{equation}
\\
Our parameter representation (19) (20) gives \\
\\
	\begin{equation}
	\begin{aligned}
	x& = s_1 - \lambda_2s_2 - \lambda_3s_3& \hspace{0.5cm} w& = s_1 + \lambda_2s_2 + \lambda_3s_3 \\
	y& = s_2 + \lambda_2s_1& \hspace{0.5cm} w& = s_2 - \lambda_2s_1 \\
	z& = s_3 + \lambda_3s_1& \hspace{0.5cm} w& = s_3 - \lambda_3s_1 \\
	\end{aligned}
	\end{equation}
\\
Conversely \\
\\
	\begin{equation}
	\begin{aligned}
	2s_1& = x + w, \hspace{0.5cm} 2s_2 = y +w, \hspace{0.5cm} 2s_3 = z + w \\
	&2s_1\lambda_2 = y - w, \hspace{0.5cm} 2s_1\lambda_3 = z - w \\
	\end{aligned}
	\end{equation}
\\
Equation (49) results in \\
\\
	\begin{equation}
	\lambda_3 = \lambda_2 + \frac{s_3 - s_2} {s_1} \\
	\end{equation}
\\
and with \\
\\
	\begin{equation}
	\begin{aligned}
	q& = s_1 + s_2 + s_3 \\
	\end{aligned}
	\end{equation}

 	\begin{equation}
	\begin{aligned}
	\hspace*{-3.25cm} \text{in} \hspace{3.75cm} \lambda_2& = \frac{1} {q} [s_2 - s_1 + s_3 \frac {s_2 - s_3} {s_1}] \\
	\end{aligned}
	\end{equation}
\\
Then \\
\\
	\begin{equation}
	\begin{aligned}
	w& = s_2 - \lambda_2s_1 \\
	w& = \frac {1} {q} [ s_1^2 + s_2^2 + s_3^3 ] \\
	\end{aligned}
	\end{equation}
\\
and with \\
\\
	\begin{equation}
	\begin{aligned}
	p = s_1s_2 + s_2s_3 + s_3s_1 \\
	\end{aligned}
	\end{equation}
\\
we find with (50) and the parameters $s_1, s_2, s_3$ \\
\\
	\begin{equation}
	\begin{aligned}
	x& = 2s_1 - w \\
	x& = \frac{1} {q} [ 2p + q ( s_1 - s_2 - s_3 ) ] \\
	\end{aligned}
	\end{equation}
\\
\\
	\begin{equation}
	\begin{aligned}
	y& = 2s_2 - w \\
	y& = \frac{1} {q} [ 2p + q ( s_2 - s_3 - s_1 ) ] \\
	\end{aligned}
	\end{equation}
\\
\\
	\begin{equation}
	\begin{aligned}
	z& = 2s_3 - w \\
	z& = \frac{1} {q} [ 2p + q ( s_3 - s_1 - s_2 ) ] \\
	\end{aligned}
	\end{equation}
\\
\\
Up to scaling and with (52) (55) we find all solutions of (48) by \\
\\
	\begin{equation}
	\begin{aligned}
	p& = s_1s_2 + s_2s_3 + s_3s_1 \\
	q& = s_1 + s_2 + s_3 \\
	x& = 2p + q ( s_1 - s_2 - s_3 ) \\
	y& = 2p + q ( s_2 - s_3 - s_1 ) \\
	z& = 2p + q ( s_3 - s_1 - s_2 ) \\
	w& = s_1^2 + s_2^2 + s_3^3 \\
	\end{aligned}
	\end{equation}
\\
Finally, introducing the six parameters \\
\\
$m_1, n_1, m_2, n_2, m_3, n_3$, all integers, through \\
\\
\centerline{$s_1 = \frac {m_1} {n_1}, \hspace{0.5cm} s_2 = \frac {m_2} {n_2}, \hspace{0.5cm} s_3 = \frac {m_3} {n_3}$,} \\
\\
we get all integer solutions of (48), up to scaling, with \\
\\
	\begin{equation}
	\begin{aligned}
	P& = m_1 n_1 m_2 n_2 n_3^2 + m_2 n_2 m_3 n_3 n_1^2 + m_3 n_3 m_1 n_1 n_2^2 \\
	Q& = m_1 n_2 n_3 + m_2 n_3 n_1 + m_3 n_1 n_2 \\
	\end{aligned}
	\end{equation}

	as follows \\

	\begin{equation}
	\begin{aligned}
	x& = 2P + Q [ m_1 n_2 n_3 - m_2 n_3 n_1 - m_3 n_1 n_2 ] \\
	y& = 2P + Q [ m_2 n_3 n_1 - m_3 n_1 n_2 - m_1 n_2 n_3 ] \\
	z& = 2P + Q [ m_3 n_1 n_2 - m_1 n_2 n_3 - m_2 n_3 n_1 ] \\
	w& = [m_1 n_2 n_3]^2 + [m_2 n_3 n_1]^2 + [m_3 n_1 n_2]^2 \\
	\end{aligned}
	\end{equation}
\\
Thus x, y, z, w are polynomials of degree six in six integer variables. \\
\\

\end{enumerate}

\bigskip

\noindent\textit{Department of Physics, University of Colorado Boulder, Boulder, CO 80309\\
Walter.Wyss@Colorado.EDU}

\end{document}